\documentclass[12pt]{article}
\usepackage{amsmath, amssymb, amsthm}

\theoremstyle{definition}

\title{Groups, periodic planes and \\
 hyperbolic buildings.}

\author{\sc Alina Vdovina}

\date{
\small {Mathematisches Institut} \\
\small {Beringstr. 1, Bonn 53115} \\
\small {e-mail: alina@math.uni-bonn.de}
}

\begin{document}

\maketitle

\abstract{
We give an elementary construction of polyhedra whose links are
connected bipartite graphs, which are not necessarily isomorphic
pairwise. We show, that the fundamental groups of some of our
polyhedra contain surface groups. In particular, we
construct polyhedra whose links are generalized $m$-gons. The
polyhedra of this type are interesting because of their universal
coverings, which are two-dimensional hyperbolic  buildings with
different links. The presentation of the results is done in the
language of combinatorial group theory.}

\bigskip

\section*{Introduction}

We will call a {\em polyhedron} a two-dimensional complex which
is obtained from several oriented $p$-gons by identification of
corresponding sides. Consider a point of the polyhedron and take
a sphere of a small radius at this point. The intersection of the
sphere with the polyhedron is a graph, which is called the {\em link}
at this point.

We consider the polyhedra such that all links of all vertices are
connected bipartite graphs.

We will say, that a polyhedron $P$ is an $(m,n)$-polyhedron,
if the girth of any link of $P$ is at least $m$ and each
face of $P$ is a polygon with at least $n$ edges.

Let $P$ be a $(m,n)$-polyhedron such that $m$ and $n$ satisfy
the inequality $mn \geq 2(m+n)$, which appear in
small cancellation theory \cite{LS}. 
The minimal solutions of the equality are $(6,3),(4,4),(3,6)$.

The universal covering of a $(m,n)$-polyhedron with
the metric introduced in \cite{BBr},p.165 is a complete
metric space of non-positive curvature in the sense
of Alexandrov and Busemann, \cite{GH}.

In chapter 2 we construct a family of $(m,n)$-polyhedra with a given even
number of sides of every face, such that the links of vertices are 
bipartite graphs, in general, nonisomorphic.


\medbreak

\noindent{\bf Definition.}
We say that sets $\mathcal{G}_1$, $\mathcal{G}_2$,..., $\mathcal{G}_k$
of connected bipartite graphs are {\em compatible}, if the number of
white vertices of all graphs of every set is equal to the same number
$n$ and there are bijections between sets of white vertices for every
pair $\mathcal{G}_1$, $\mathcal{G}_j$, $j=2,...,k$, preserving the degrees
of vertices. 

\medbreak

\noindent {\bf Theorem 1.}
Let $\mathcal{G}_1$, $\mathcal{G}_2$,..., $\mathcal{G}_k$ be compatible
sets of connected bipartite graphs and $k\ge1$. Then there exists a
family of finite polyhedra with $2k$-gonal faces and links at vertices
isomorphic to the graphs from $\mathcal{G}_1$,
$\mathcal{G}_2$,..., $\mathcal{G}_k$.

\bigbreak

Let $X$ be a universal covering of a $(m,n)$-polyhedron $P$.
Then, by a result of Gromov \cite{Gr}, p.119, the fundamental
group $G$ of $P$ is hyperbolic iff $X$ does not contain a flat
Euclidean plane. A natural question is the following:
 whether $G$ contains
 $\mathbb{Z} \oplus \mathbb{Z}$ if $X$ contains a flat plane.

This question is a particular case of more general problem,
formulated by Ballmann and Brin for polygonal
$(m,n)$-complexes \cite{BBr},p.166. They solved it
for all $(3,6)$-complexes \cite{BBr} and for $(6,3)$-complexes
in the case of euclidean buildings \cite{BB}.

In chapter 3 we show that fundamental group of our polyhedron contains
$\mathbb{Z} \oplus \mathbb{Z}$  if its universal covering
contains a flat plane.

In chapter 4 we apply the construction from Chapter 2
to prove existence of periodic planes in some hyperbolic buildings.

\bigbreak

\noindent{\bf Definition.}
A {\em generalized $m$-gon} is a graph which is a union of subgraphs,
called apartments, such that:
\begin{itemize}
\item[1.] Every apartment is a cycle composed of $2m$ edges
for some fixed $m$.
\item[2.] For any two edges there is an apartment containing both of them.
\item[3.] If two apartments have a common edge, then there is
an isomorphism between them fixing their intersection pointwise.
\end{itemize}

\bigbreak

\noindent
{\bf Definition.} Let $\mathcal{P}(p,m)$ be a tessellation of the
hyperbolic plane by  regular polygons with $p$ sides,
with angles $\pi/m$ in each vertex where $m$ is an integer.
A {\em hyperbolic building} is a polygonal complex $X$,
which can be expressed as the union of subcomplexes called apartments 
such that:

\begin{itemize}
\item[1.] Every apartment is isomorphic to $\mathcal{P}(p,m)$.
\item[2.] For any two polygons of $X$, there is an apartment
containing both of them.
\item[3.] For any two apartments $A_1, A_2 \in X$ containing
the same polygon, there exists an isomorphism $ A_1 \to A_2$
fixing $A_1 \cap A_2$.
\end{itemize}

\bigbreak

Let $C_p$ be a polyhedron whose faces are $p$-gons
and whose links are generalized $m$-gons with $mp>2m+p$.
We equip every face of $C_p$ with the hyperbolic
metric such that all sides of the polygons are geodesics
and all angles are $\pi/m$. Then the universal covering of such
polyhedron is a hyperbolic building, see \cite{Paulin}.

In the case $p=3$, $m=3$, i.e. $C_p$ is a simplex, we can give
a euclidean metric to every face. In this metric
all sides of the triangles
are geodesics. The universal coverings of these polyhedra
with the euclidean metric are euclidean buildings, see \cite{BB}, \cite{Ba}.
The metric characterization of euclidean buildings can be found in \cite{CL}.




It is known (cf. \cite{Ronan}), that 
a {\em generalized m-gon} is a connected, bipartite graph of diameter $m$
and girth $2m$, in which each vertex lies
on at least two edges. A graph is {\em bipartite} if its set of vertices
can be partitioned into two disjoint subsets such that no
two vertices in the same subset lie on a common edge.
The vertices of the one subset we will call black vertices
and the vertices of the other subset the white ones.
The {\em diameter} is the maximum distance between two vertices
and the {\em girth}  is the length of a shortest circuit.

Let $G$ be a connected bipartite graph on $q+r$ vertices,
$q$ black vertices and $r$ white ones. Let $\mathcal{A}$ and
$\mathcal{B}$ be two alphabets on $q$ and $r$ letters respectively,
$\mathcal{A}=\{x_1, x_2, \ldots, x_q\}$ and $\mathcal{B}=\{y_1, y_2, \ldots,
y_r\}$.
We mark every black vertex with an element from
$\mathcal{A}$ and every white vertex with an element from $\mathcal{B}$.

We will  also define an incidence tableau of such a graph in the following
way: the first element of each line is a white vertex
and all other elements are black vertices incident to this white vertex.
Different lines correspond to different white vertices.

\bigbreak

\noindent {\bf Example.} For the smallest complete bipartite graph
(generalized 2-gon) with four vertices,
two black $x_1, x_2$  and two white $y_1, y_2$,
the incidence matrix will be

$y_1$:  $x_1$ $x_2$

$y_2$:  $x_1$ $x_2$

\bigbreak

Define a graph $G^{\prime}$, called the {\em dual} to $G$,
in the following way: the graph $G^{\prime}$ can be obtained from $G$
by changing black vertices to white and vice-versa.
(It is easy to see that $G$ and $G^{\prime}$ are isomorphic
as ordinary graphs.)

 
Similarly, we can define the incidence tableau
by "inverse" way, as correspondence of black vertices
to white ones. For the graph from the example we have

$x_1$:  $y_1$ $y_2$

$x_2$:  $y_1$ $y_2$


\section{Polygonal presentation.}

\noindent
{\bf Definition.} Suppose we have $n$ disjoint connected bipartite graphs
 $G_1, G_2, \ldots G_n$.
Let $P_i$ and $Q_i$ be the sets of black and white vertices respectively in
$G_i$, $i=1,...,n$; let $P=\cup P_i, Q=\cup Q_i$, $P_i \cap P_j = \emptyset$
 $Q_i \cap Q_j = \emptyset$
for $i \neq j$ and
let $\lambda$ be a  bijection $\lambda: P\to Q$.

A set $\mathcal{K}$ of $k$-tuples $(x_1,x_2, \ldots, x_k)$, $x_i \in P$,
will be called a {\em polygonal presentation} over $P$ compatible
with $\lambda$ if

\begin{itemize}

\item[(1)] $(x_1,x_2,x_3, \ldots ,x_k) \in \mathcal{K}$ implies that
   $(x_2,x_3,\ldots,x_k,x_1) \in \mathcal{K}$;

\item[(2)] given $x_1,x_2 \in P$, then $(x_1,x_2,x_3, \ldots,x_k) \in \mathcal{K}$
for some $x_3,\ldots,x_k$ if and only if $x_2$ and $\lambda(x_1)$
are incident in some $G_i$;

\item[(3)] given $x_1,x_2 \in P$, then  $(x_1,x_2,x_3, \ldots ,x_k) \in \mathcal{K}$
    for at most one $x_3 \in P$.

\end{itemize}

If there exists such $\mathcal{K}$, we will call $\lambda$ a {\em basic bijection}.

Polygonal presentations for $n=1$, $k=3$ were listed in $\cite{Cart}$
with the incidence graph of the finite projective plane of order two or three as the graph $G_1$.
 They were called  triangle
presentations in that paper, because the case $k=3$ was considered.
Polygonal presentations for $k > 3$ were considered in \cite{V}.

\section{Construction of polyhedra.}

We can associate  a polyhedron $K$ on $n$ vertices with
each polygonal presentation $\mathcal{K}$ as follows:
for every cyclic $k$-tuple $(x_1,x_2,x_3,\ldots,x_k)$ from
the definition
we take an oriented $k$-gon on the boundary of which
the word $x_1 x_2 x_3\ldots x_k$ is written. To obtain
the polyhedron we identify the corresponding sides of our
polygons, respecting orientation.
We will say that the
polyhedron $K$ {\em corresponds} to the polygonal
presentation $\mathcal{K}$.
\medskip

\noindent
{\bf Lemma\cite{V}} A polyhedron $K$ which corresponds to
a polygonal presentation $\mathcal{K}$ has
  graphs $G_1, G_2, \ldots, G_n$ as the links.

\medskip

\noindent
{\bf Remark.} Consider a  polygonal
presentation $\mathcal{K}$. Let $s_i$ be the number of vertices
of the graph $G_i$ and $t_i$ be the number of edges of $G_i$,
$i=1,...,n$.
If the polyhedron $K$  corresponds to the polygonal
presentation $\mathcal{K}$, then $K$ has $n$ vertices
(the number of vertices of $K$ is equal to the number of graphs),
$k \sum_{i=1}^n s_i$ edges and $\sum_{i=1}^n t_i$ 
faces, all faces are polygons with $k$ sides.

\medskip
 
\noindent
{\bf Proof of Theorem 1.}
By the Lemma, to construct the polyhedron with given links,
it is sufficient to construct a corresponding polygonal presentation.

By the definition of compatible sets of bipartite graphs,
there are bijections $\alpha_j, j=2,...,k,$ from the set
of white vertices of  $\mathcal{G}_1$
to the set of white vertices $\mathcal{G}_j$ preserving the degrees
of white vertices.

We mark white vertices of $\mathcal{G}_i, i=1,...,k,$ by letters
of an alphabet $\mathcal{A}_i, i=1, \ldots, 2k$,
 $\mathcal{A}_i=\{ x_1^i, \ldots, x_q^i\}$, such that the bijections
$\alpha_j, j=2,...,k,$ are induced by the indexes of letters
, i.e. $\alpha_j(x_m^1)=x_m^j$,
 $ j=2,...,k,$ 
We mark black vertices of $\mathcal{G}_i, i=1,...,k$ by letters
of an alphabet $\mathcal{B}_i=\{y_1^i, y_2^i, \ldots,
y_r^i\}$. So, every edge of $\mathcal{G}_i, i=1,...,k$  can be presented in
a form $(x_m^iy_l^i)$, $m=1,...,n, l=1,...,r$.

Having such a set of bijections $\alpha_j, j=2,...,k$ of white
vertices we can choose bijections $\beta_j, j=2,...,k$
of the set of edges of $\mathcal{G}_1$  
to the set of edges of $\mathcal{G}_j$ which preserves
$\alpha_j, j=2,...,k$. 
Let $\beta_j(x^1_m y^1_l)=x_m^j y_{l_j}^j$,
then we take the cyclic word $(x^1_m, y^1_{l_1},...,x_m^j, y_{l_j}^j,...,
x_m^k ,y_{l_k}^k)$, $m=1,...,n, l_t=1,...,r, t=1,...,k$,
to the set $\mathcal{P}$.
We will prove, that $\mathcal{P}$ is a polygonal presentation.

For this we need $k$ more families of graphs $\mathcal{H}_i, i=1,...,k$,
such that every graph $H$ is contained in $\mathcal{H}_i$ 
if and only if it is dual for a graph $G$
from $\mathcal{G}_i$, $i=1,...,k$, and every graph from $\mathcal{G}_i$,
 $i=1,...,k$, has its dual in $\mathcal{H}_i, i=1,...,k$.

We mark black vertices of $\mathcal{H}_i, i=1,...,k$ by letters
of an alphabet $\mathcal{A}_i, i=1, \ldots, 2k$,
 $\mathcal{A}_i=\{ x_1^i, \ldots, x_q^i\}$ and
we mark white vertices of $\mathcal{H}_i, i=1,...,k$ by letters
of an alphabet $\mathcal{B}_i=\{y_1^i, y_2^i, \ldots,
y_r^i\}$.

The bijection $\lambda$ from the set $P$ of all black vertices
of graphs from $\mathcal{G}_j$, $j=1,...,k$ and $\mathcal{H}_j, j=1,...,k$
to the set $Q$ of all white ones is defined by labels:
$\lambda(x_m^j)=x_m^j$ and $\lambda(y_l^j)=y_l^i$.

It is necessary to check all axioms of the polygonal presentation
for an arbitrary $2k$-tuple $(x^1_m, y^1_{l_1},...,x_m^j, y_{l_j}^j,...,
x_m^k ,y_{l_k}^k)$, $m=1,...,n, l_t=1,...,r, t=1,...,k$:

\begin{itemize}

\item[(1)]
In our construction we take all cyclic permutations
of each $2k$-tuple.

\item[(2)]
Let's consider an arbitrary $2k$-tuple from $\mathcal{P}$.
There are two possibilities: when the tuple starts with an element
from  $\mathcal{A}_j$ or $\mathcal{B}_j$, namely
 $(x_m^j, y_{l_j}^j,...,
x_m^k ,y_{l_k}^k,...,x^{j-1}_m, y^{j-1}_{l_1})$ or 
 
$(y_{l_j}^j,x_m^{j+1},..,
x_m^k ,y_{l_k}^k...x^{j-1}_m, y^{j-1}_{l_1})$.
In the case of 

$(x_m^j, y_{l_j}^j,...,
x_m^k ,y_{l_k}^k...x^{j-1}_m, y^{j-1}_{l_1})$ we have, that
$\lambda(x_m^j)$ is a white vertex $x_m^j$ of some graph $G$
 from $\mathcal{G}_j$ and  $y_{l_j}^j$ is a black vertex of $G$
and  $x_m^j$ and $y_{l_j}^j$ are incident in $G$ by the construction
of $\mathcal{P}$.

In the case of $(y_{l_j}^j,x_m^{j+1},..,
x_m^k ,y_{l_k}^k...x^1_m, y^1_{l_1})$
$\lambda(y_{l_j}^j)$ is a white vertex $y_{l_j}^j$ of some graph $H$
 from $\mathcal{H}_j$ and  $x_m^{j+1}$ is a black vertex of $H$
and  $x_m^{j+1}$ and $y_{l_j}^j$ are incident in $H$ by construction
of $\mathcal{P}$.

\item[(3)]
 Since $\alpha_j, j=2,...,k$ and $\beta_j, j=2,...,k$
are bijections, then there is a unique word in  $\mathcal{P}$,
which contains given subword of length two. This proves
the property (3) of the polygonal presentation.

\end{itemize}

The Theorem is proved.

\bigbreak

\noindent We will denote polyhedra constructed in the proof of the Theorem 1
$P_k$. 


\section{$\mathbb{Z} \oplus \mathbb{Z}$ as a subgroup of the fundamental group
of $P_2$.}

 In this chapter we  give
an euclidean metric to every face of $P_2$. In this metric
all sides of the squares are geodesics and angles are $\pi/2$.

{\bf Definition.}
A plane in the universal covering of $P_2$ is {\em periodic},
if it is stabilized  by a $\mathbb{Z} \oplus \mathbb{Z}$
subgroup of the fundamental group of $P_2$.


\bigbreak

\noindent
{\bf Theorem 2.} If the universal covering $U$ of $P_2$ contains
a flat plane $F$, then $U$ contains a periodic plane and $\pi (P_2)$
contains a subgroup isomorphic to
 $\mathbb{Z} \oplus \mathbb{Z}$, where 
$\pi (P_2)$ is the fundamental group of $P_2$.

\bigbreak
\begin{picture}(310,90)
\put(5,5){\line(1,0){5}}
\put(5,45){\line(1,0){5}}
\put(5,85){\line(1,0){5}}
\put(10,5){\vector(0,1){23}}\put(10,28){\line(0,1){17}}\put(12,27){$v_{58}$}
\put(10,85){\vector(0,-1){23}}\put(10,62){\line(0,-1){17}}\put(12,64){$v_{25}$}
\put(50,5){\vector(-1,0){23}}\put(27,5){\line(-1,0){17}}\put(26,8){$u_1$}
\put(10,45){\vector(1,0){23}}\put(33,45){\line(1,0){17}}\put(29,48){$x_1$}
\put(50,85){\vector(-1,0){23}}\put(27,85){\line(-1,0){17}}\put(23,78){$u_1$}
\put(50,45){\vector(0,-1){23}}\put(50,22){\line(0,-1){17}}\put(39,21){$y_2$}
\put(50,45){\vector(0,1){23}}\put(50,68){\line(0,1){17}}\put(39,65){$y_1$}
\put(50,45){\circle*{3}}
\put(50,5){\vector(1,0){23}}\put(73,5){\line(1,0){17}}\put(68,9){$u_2$}
\put(90,45){\vector(-1,0){23}}\put(67,45){\line(-1,0){17}}\put(66,49){$x_2$}
\put(50,85){\vector(1,0){23}}\put(73,85){\line(1,0){17}}\put(70,78){$u_2$}
\put(90,5){\vector(0,1){23}}\put(90,28){\line(0,1){17}}\put(75,27){$v_{88}$}
\put(90,85){\vector(0,-1){23}}\put(90,62){\line(0,-1){17}}\put(75,60){$v_{73}$}
\put(90,5){\line(1,0){60}}
\put(90,45){\vector(1,0){23}}\put(113,45){\line(1,0){17}}\put(109,49){$x_3$}
\put(90,85){\line(1,0){30}}
\put(160,5){\line(1,0){50}}
\put(150,45){\line(1,0){60}}
\put(130,85){\line(1,0){80}}
\put(170,5){\line(0,1){80}}
\put(135,44){...}
\put(210,5){\vector(0,1){23}}\put(210,28){\line(0,1){17}}\put(212,27){$v_{58}$}
\put(210,85){\vector(0,-1){23}}\put(210,62){\line(0,-1){17}}\put(212,64){$v_{25}$}
\put(250,5){\vector(-1,0){23}}\put(227,5){\line(-1,0){17}}\put(226,8){$u_1$}
\put(210,45){\vector(1,0){23}}\put(233,45){\line(1,0){17}}\put(229,48){$x_1$}
\put(250,85){\vector(-1,0){23}}\put(227,85){\line(-1,0){17}}\put(223,78){$u_1$}
\put(250,45){\vector(0,-1){23}}\put(250,22){\line(0,-1){17}}\put(239,21){$y_2$}
\put(250,45){\vector(0,1){23}}\put(250,68){\line(0,1){17}}\put(239,65){$y_1$}
\put(250,45){\circle*{3}}
\put(250,5){\vector(1,0){23}}\put(273,5){\line(1,0){17}}\put(268,9){$u_2$}
\put(290,45){\vector(-1,0){23}}\put(267,45){\line(-1,0){17}}\put(266,49){$x_2$}
\put(250,85){\vector(1,0){23}}\put(273,85){\line(1,0){17}}\put(270,78){$u_2$}
\put(290,5){\vector(0,1){23}}\put(290,28){\line(0,1){17}}\put(275,27){$v_{88}$}
\put(290,85){\vector(0,-1){23}}\put(290,62){\line(0,-1){17}}\put(275,60){$v_{73}$}
\put(290,5){\line(1,0){20}}
\put(290,45){\line(1,0){15}}
\put(290,85){\line(1,0){10}}
\end{picture}

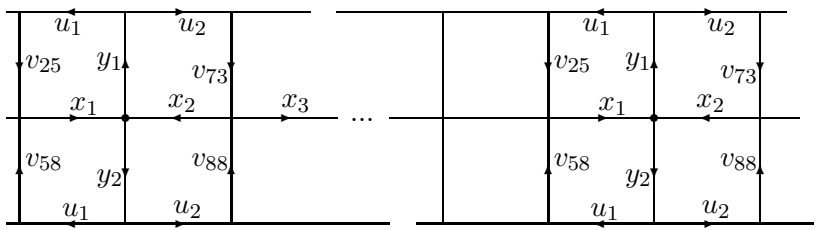
\centerline{\bf Figure 1.}

\medskip

\noindent{\bf Proof.}
All words of the polygonal presentation for $P_2$
have form $x_iy_ju_iv_l$. So, in the flat plane $F$ there are
$x,y,u,v$-lines. Let's consider an infinite strip $S$
consisting of $x$-line $L$ and all squares, which contain
common edges from $L$ (fig.1). Let $v$ be an internal vertex
of $S$ with entering edges $x_1, x_2$ and leaving edges $y_1, y_2$.  
Because of the finiteness of the polyhedron, there is a finite
number of $x$'s and $y$'s. So, in our infinite strip $S$
there is another internal vertex $w$  with entering edges 
$x_1, x_2$ and leaving edges $y_1, y_2$. Let the distance
between $v$ and $w$ be $s$(the distance is
the number of $x$-edges). Let's consider the rectangle $R$
which consists of a part
of $L$ between $v$ and $w$ and all squares, which contain
common edges from $L$ (fig.1). Consider the $x$-word $U$ of length $s$   
$x_2x_3...x_1$ 
written
between $v$ and $w$. We define the word $W$ to be obtained from $U$
by replacing $x$ by $u$. Then, on the boundary of $R$ we can read
a word $y_2^{-1}y_1W^{-1}y_1^{-1}y_2W$. We consider a plane
$F^{\prime}$ in $U$ tesselated by $R$. This tesselation
exists, because the sum of all angles at any vertex of it
is equal to $2 \pi$ (the degree of every vertex is four and each
angle is $\pi/2$).
The subgroup of  $\pi (P_2)$, generated by $a=y_2^{-1}y_1$ and $b=W^{-1}$
is isomorphic to 
$\mathbb{Z} \oplus \mathbb{Z}$ and acts on $F^{\prime}$ uniformly.

The Theorem 2 is proved.

\medbreak

Let's note, that polyhedra of type $P_2$ are particular cases
of polygonal complexes $(4,4)$, which were considered in \cite{BBr}.



\section{Periodic planes in hyperbolic buildings.}


In this chapter we consider polyhedra whose faces are $2k$-gons
and whose links are generalized $m$-gons with $km>m+k$.
We equip every face of the polyhedra with the hyperbolic
metric such that all sides of the polygons are geodesics
and all angles are $\pi/m$. Then the universal covering of such
polyhedra are hyperbolic buildings, see \cite{Paulin}.

\medbreak

In the hyperbolic case the notion of periodicity should
be extended to the action of surface groups of genus $g \geq 2$.

{\bf Definition.}
We will say, that a tesselated plane is {\em periodic}, if there exists
a genus $g \geq 1$ sufrace group acting on it uniformly.

It is natural question whether fundamental group of a hyperbolic
polyhedron always contains a surface group, if its universal
covering contains a hyperbolic plane.
 In this chapter we will prove that
this is true for some hyperbolic buildings.
For the proofs we will need some notations concerning
Wicks forms.

\medbreak

\noindent{\bf Definition }
 An {\it oriented Wicks form\/} is a cyclic word $w= w_1w_2\dots w_{2l}$
 (a cyclic word is the orbit of a linear word under cyclic permutations)
 in some alphabet $a_1^{\pm 1},a_2^{\pm 1},\dots$ of letters
 $a_1,a_2,\dots$ and their inverses $a_1^{-1},a_2^{-1},\dots$ such that
\begin{itemize}
\item[(i)] if $a_i^\epsilon$ appears in $w$ (for $\epsilon\in\{\pm 1\}$)
 then $a_i^{-\epsilon}$ appears exactly once in $w$,
\item[(ii)] the word $w$ contains no cyclic factor (subword of
 cyclically consecutive letters in $w$) of the form $a_i a_i^{-1}$ or
 $a_i^{-1}a_i$ (no cancellation),
\item[(iii)] if $a_i^\epsilon a_j^\delta$ is a cyclic factor of $w$ then
 $a_j^{-\delta}a_i^{-\epsilon}$ is not a cyclic factor of $w$ 
(substitutions of the form
 $a_i^\epsilon a_j^\delta\longmapsto x,
 \quad a_j^{-\delta}a_i^{-\epsilon}\longmapsto x^{-1}$ are impossible).
\end{itemize}

 An oriented Wicks form $w=w_1w_2\dots$ in an alphabet $A$ is
 {\em isomorphic\/} to $w'=w'_1w'_2$ in an alphabet $A'$ if
 there exists a bijection $\varphi:A\longrightarrow A'$ with
 $\varphi(a^{-1})=\varphi(a)^{-1}$ such that $w'$ and
 $\varphi(w)=\varphi(w_1)\varphi(w_2)\dots$ define the
 same cyclic word.

 An oriented Wicks form $w$ is an element of the commutator subgroup
 when considered as an element in the free group $G$ generated by
 $a_1,a_2,\dots$. We define the {\em algebraic genus\/} $g_a(w)$ of
 $w$ as the least positive integer $g_a$ such that $w$ is a product
 of $g_a$ commutators in $G$.

 The {\em topological genus\/} $g_t(w)$ of an oriented Wicks
 form $w=w_1\dots w_{2e-1}w_{2e}$ is defined as the topological
 genus of the oriented compact connected surface obtained by
 labeling and orienting the edges of a $2e-$gone (which we
 consider as a subset of the oriented plane) according to
 $w$ and by identifying the edges in the obvious way.

\bigbreak

\noindent{\bf Proposition }
{\sl The algebraic and the topological genus of an oriented Wicks
 form coincide (cf. \cite{C},\cite{CE}).}

\bigbreak

 We define the {\em genus\/} $g(w)$ of an oriented
 Wicks form $w$ by $g(w)=g_a(w)=g_t(w)$.

 Consider the oriented compact surface $S$ associated to an oriented 
 Wicks form $w=w_1\dots w_{2e}$. This surface carries an immerged graph
 $\Gamma\subset S$ such that $S\setminus \Gamma$ is an open polygon
 with $2e$ sides (and hence connected and simply connected).
 Moreover, conditions (ii) and (iii) on Wicks form imply that $\Gamma$ 
 contains no vertices of degree $1$ or $2$ (or equivalently that the
 dual graph of $\Gamma\subset S$ contains no faces which are $1-$gones
 or $2-$gones). This construction works also
 in the opposite direction: Given a graph $\Gamma\subset S$
 with $e$ edges on an oriented compact connected surface $S$ of genus $g$
 such that $S\setminus \Gamma$ is connected and simply connected, we get
 an oriented Wicks form of genus $g$ and length $2e$ by labeling and 
 orienting the edges of $\Gamma$ and by cutting $S$ open along the graph
 $\Gamma$. The associated oriented Wicks form is defined as the word
 which appears in this way on the boundary of the resulting polygon
 with $2e$ sides. We identify henceforth oriented Wicks
 forms with the associated immerged graphs $\Gamma\subset S$,
 speaking of vertices and edges of oriented Wicks form.

 The formula for the Euler characteristic
 $$\chi(S)=2-2g=v-e+1$$
 (where $v$ denotes the number of vertices and $e$ the number
 of edges in $\Gamma\subset S$) shows that
 an oriented Wicks form of genus $g$ has at least length $4g$
 (the associated graph has then a unique vertex of degree $4g$
 and $2g$ edges) and at most length $6(2g-1)$ (the associated
 graph has then $2(2g-1)$ vertices of degree three and
 $3(2g-1)$ edges).

\medbreak

\noindent {\bf Definition.} We will say that word $W$
is obtained from a word $U$ by a non-cancelling substitution $\Phi$,
if we substitute every letter of $U$ such that there is no
cancellations between $\Phi(y_i)$ and $\Phi(y_j)$
whenever $y_iy_j$ is a subword of $U$.

\medbreak

Obviously,  from a disc with the word
$U$ on its boundary we can get a genus $g$ surface by identification
of sides (with the same labels respecting orientation) if and only
if a word $U$ is obtained from a Wicks form of genus
$g$ by non-cancelling substitution.

\medbreak

\noindent {\bf Theorem 3.} If $I_{2k}, k \geq 3$ is a right-angled hyperbolic
building which apartments are hyperbolic planes tesselated by
polygons with $2k$ sides, then $I_{2k}$ contains
a periodic plane under the action of genus $g=2k-4$ surface group.

\noindent
\begin{picture}(300,300)
\put(150,150){\circle*{3}}\put(152,152){$v$}
\put(150,150){\vector(0,-1){30}}\put(150,120){\line(0,-1){90}}\put(155,113){$y_j^1$}
\put(150,40){\vector(-1,-1){10}}\put(140,30){\line(-1,-1){10}}\put(135,40){$x_i^2$}
\put(130,10){\line(0,1){20}}
\put(130,20){\vector(-1,0){40}}\put(90,20){\line(-1,0){30}}\put(85,6){$y_j^2$}
\put(70,150){\vector(1,0){60}}\put(130,150){\line(1,0){20}}\put(110,155){$x_i^1$}
\put(30,100){\vector(1,1){30}}\put(60,130){\line(1,1){20}}\put(60,120){$y_j^k$}
\put(30,60){\vector(0,1){25}}\put(30,85){\line(0,1){15}}\put(33,80){$x_i^k$}
\put(20,100){\line(1,0){20}}
\put(45,35){\circle*{2}}\put(40,40){\circle*{2}}\put(35,45){\circle*{2}}
\put(150,40){\vector(1,-1){10}}\put(160,30){\line(1,-1){10}}\put(162,30){$x_s^2$}
\put(170,10){\line(0,1){15}}
\put(170,20){\vector(1,0){40}}\put(210,20){\line(1,0){30}}\put(215,6){$y_j^2$}
\put(230,150){\vector(-1,0){60}}\put(170,150){\line(-1,0){20}}\put(175,155){$x_s^1$}
\put(270,100){\vector(-1,1){30}}\put(240,130){\line(-1,1){20}}\put(250,130){$y_j^k$}
\put(270,60){\vector(0,1){25}}\put(270,85){\line(0,1){15}}\put(273,80){$x_s^k$}
\put(250,30){\circle*{2}}\put(255,35){\circle*{2}}\put(260,40){\circle*{2}}
\put(260,90){\line(1,1){20}}
\put(150,150){\vector(0,1){50}}\put(150,200){\line(0,1){70}}\put(155,190){$y_t^1$}
\put(150,260){\vector(1,1){10}}\put(160,270){\line(1,1){10}}\put(160,260){$x_s^2$}
\put(165,285){\line(1,-1){10}}
\put(170,280){\vector(1,0){40}}\put(210,280){\line(1,0){30}}\put(215,268){$y_t^2$}
\put(270,200){\vector(-1,-1){30}}\put(240,170){\line(-1,-1){20}}\put(250,170){$y_t^k$}
\put(270,240){\vector(0,-1){25}}\put(270,215){\line(0,-1){15}}\put(273,220){$x_s^k$}
\put(255,270){\circle*{2}}\put(260,265){\circle*{2}}\put(265,260){\circle*{2}}
\put(150,260){\vector(-1,1){10}}\put(140,270){\line(-1,1){10}}\put(125,260){$x_i^2$}
\put(265,205){\line(1,-1){10}}
\put(125,275){\line(1,1){10}}
\put(25,195){\line(1,1){10}}
\put(130,280){\vector(-1,0){40}}\put(90,280){\line(-1,0){30}}\put(90,268){$y_t^2$}
\put(30,200){\vector(1,-1){30}}\put(60,170){\line(1,-1){20}}\put(60,175){$y_t^k$}
\put(30,240){\vector(0,-1){25}}\put(30,215){\line(0,-1){15}}\put(33,220){$x_i^k$}
\put(35,255){\circle*{2}}\put(40,260){\circle*{2}}\put(45,265){\circle*{2}}
\put(260,0){\bf{Figure 2.}}
\end{picture}

\bigbreak

\noindent {\bf Proof.} Since a right-angled hyperbolic building with
given local data is unique ( \cite{B1} \cite{B2}),
then it can be obtained as a universal covering of a polyhedron,
described in  \cite{V}.
As $G$ it is necessary to take a complete bipartite graph.
This is also a particular case of the construction from
the chapter 1, when every set of graphs consists from the
same graph $G$, which is a complete bipartite
graph. Let's consider any apartment in $I_{2k}, k \geq 3$.
It is a hyperbolic plane tesselated by regular right-angled polygons with
$2k$ sides, such that all vertices of the tesselation
have degree four. Consider a vertex $v$ such that edges $x^1_i$
and $x^1_s$ are entering $v$ and $y^1_t$ and $y^1_j$ leaving $v$
(fig.2).There are exactly four polygons containing $v$:

$x_i^1y_j^1x_i^2y_j^2...x_i^ky_j^k$,
$x_i^1y_t^1x_i^2y_t^2...x_i^ky_t^k$,

$x_s^1y_j^1x_s^2y_j^2...x_s^ky_j^k$,
$x_s^1y_t^1x_s^2y_t^2...x_s^ky_t^k$.

Consider the word 
$$W=x_i^2y_j^2...x_i^ky_j^k(x_i^2y_t^2...x_i^ky_t^k)^{-1}
x_s^2y_t^2...x_s^ky_t^k(x_s^2y_j^2...x_s^ky_j^k)^{-1},
$$ written on the boundary of the region $D$,
consisting of polygons

$x_i^1y_j^1x_i^2y_j^2...x_i^ky_j^k$,
$x_i^1y_t^1x_i^2y_t^2...x_i^ky_t^k$,

$x_s^1y_j^1x_s^2y_j^2...x_s^ky_j^k$,
$x_s^1y_t^1x_s^2y_t^2...x_s^ky_t^k$(fig.2).
The word $W$ can be obtained from a Wicks form 
$$U=a_i^2b_j^2...a_i^kb_j^k(a_i^2b_t^2...a_i^k)^{-1}
b_t^2...a_s^k(b_j^2...a_s^kb_j^k)^{-1},$$

by a non-cancelling substitution, defined as follows:
$a_i^n=x_i^n, n=3,...,k, a_s^n=x_s^n, n=3,...,k,
b_i^n=y_i^n, n=2,...,k-1, b_t^n=y_t^n, n=2,...,k-1,
b_j^k=y_j^k (y_t^k)^{-1},
a_i^2=(x_s^2)^{-1}x_i^2$. The graph of $U$ is on fig.3.

\begin{picture}(300,300)
\put(150,150){\circle*{5}}
\put(150,150){\vector(1,1){30}}\put(180,180){\line(1,1){20}}\put(170,185){$a_s^k$}
\put(150,150){\vector(-1,1){30}}\put(120,180){\line(-1,1){20}}\put(120,185){$a_i^k$}
\put(150,200){\oval(100,20)[t]}
\put(150,210){\circle*{2}}\put(150,230){\circle{40}}\put(145,235){$b_j^k$}
\put(100,100){\vector(1,1){30}}\put(130,130){\line(1,1){20}}\put(130,120){$b_t^2$}
\put(200,100){\vector(-1,1){30}}\put(170,130){\line(-1,1){20}}\put(160,120){$b_j^2$}
\put(150,100){\oval(100,20)[b]}
\put(150,90){\circle*{2}}\put(150,70){\circle{40}}\put(145,55){$a_i^2$}
\put(150,150){\line(5,4){100}}\put(150,150){\line(5,1){100}}\put(250,200){\oval(30,60)[r]}
\put(150,150){\line(5,-4){100}}\put(150,150){\line(2,-1){100}}\put(250,85){\oval(25,30)[r]}
\multiput(250,145)(0,-7){3}{\circle*{2}}
\put(150,150){\line(-5,4){100}}\put(150,150){\line(-5,1){100}}\put(50,200){\oval(30,60)[l]}
\put(150,150){\line(-5,-4){100}}\put(150,150){\line(-2,-1){100}}\put(50,85){\oval(25,30)[l]}
\multiput(50,145)(0,-7){3}{\circle*{2}}
\put(120,20){\bf{Figure 3.}}
\end{picture}

We can tesselate the hyperbolic plane by the region $D$, according
to the word $W$, since $W$ is quadratic and the sum of all angles
at every vertex of the tesselation is $2 \pi$.

The graph of $U$ has $4k-6$ edges and $3$ vertices.
The formula for the Euler characteristic
gives $3-(4k-6)+1=2-2g$, where $g$ is the genus of $U$. So 
$g=2k-4$.

The genus of the word $U$ and therefore 
of $W$ is $2k-4$ and the hyperbolic plane
tesselated by $D$ is a periodic plane under the action of genus 
$g=2k-4$ surface group.

The Theorem is proved.


\bigbreak

\noindent {\bf Theorem 4.}
Let $I_{m,4}$ be a hyperbolic building, obtained from
the hyperbolic polyhedron $P_4$ with a generalized $m$-gon
as a link, $m=3,4,6,8$, then $I_{m,4}$ contains
a periodic plane under the action of genus $g=m-1$ surface group.

\smallbreak

\begin{picture}(300,170)
\put(30,0){\begin{picture}(230,170)%
\put(30,30){\vector(1,0){30}}\put(60,30){\line(1,0){20}}\put(57,33){$b_m$}
\put(80,30){\vector(1,1){20}}\put(100,50){\line(1,1){20}}\put(98,40){$c_m$}
\put(120,70){\vector(-1,0){30}}\put(90,70){\line(-1,0){20}}\put(77,59){$d_m$}
\put(70,70){\vector(-1,-1){20}}\put(50,50){\line(-1,-1){20}}\put(35,53){$a_m$}
\put(70,70){\vector(-1,1){20}}\put(50,90){\line(-1,1){20}}\put(38,85){$a_1$}
\put(30,110){\vector(1,0){30}}\put(60,110){\line(1,0){20}}\put(58,115){$b_m$}
\put(80,110){\vector(1,-1){20}}\put(100,90){\line(1,-1){20}}\put(90,85){$c_1$}
\put(120,70){\vector(1,1){20}}\put(140,90){\line(1,1){20}}\put(130,95){$d_1$}
\put(160,110){\vector(-1,1){20}}\put(140,130){\line(-1,1){20}}\put(140,135){$a_1$}
\put(120,150){\vector(-1,-1){20}}\put(100,130){\line(-1,-1){20}}\put(90,135){$b_1$}
\put(160,110){\vector(1,0){30}}\put(190,110){\line(1,0){20}}\put(180,115){$a_2$}
\put(210,110){\vector(-1,-1){20}}\put(190,90){\line(-1,-1){20}}\put(195,85){$b_1$}
\put(170,70){\vector(-1,0){30}}\put(140,70){\line(-1,0){20}}\put(155,75){$c_2$}
\put(120,70){\vector(1,-1){20}}\put(140,50){\line(1,-1){20}}\put(140,55){$d_2$}
\put(160,30){\vector(1,0){30}}\put(190,30){\line(1,0){20}}\put(177,35){$a_2$}
\put(210,30){\vector(-1,1){20}}\put(190,50){\line(-1,1){20}}\put(190,55){$b_2$}
\multiput(118,30)(5,0){3}{\circle*{2}}
\put(120,70){\circle*{4}}
\put(100,5){\bf{Figure 4.}}
\end{picture}}
\end{picture}

\smallbreak

\noindent {\bf Proof.} Let's consider any apartment in $I_{m,4}$.
It is a hyperbolic plane tesselated by polygons with
$4$ sides, such that all vertices of the tesselation
have degree $2m$. Consider a vertex $v$ such that edges $c_1,...,c_m$
 are entering $v$, $d_1,...,d_m$ leaving $v$ and entering and leaving
edges alternate (fig.4). There are exactly $2m$ polygons containing $v$:
$a_1b_1c_1d_1$, $a_1b_2c_1d_2$, $a_2b_2c_2d_2$,...,
$a_m b_m c_m d_m$,  $a_m b_1 c_m d_m$. Consider the region $D$
consisting of those polygons and the word 
$$W=a_1b_1b_2^{-1}a_1^{-1}a_2...a_mb_mb_1^{-1}a_m^{-1},$$ written on the
boundary of $D$. 

We can tesselate the hyperbolic plane by the region $D$, according
to the word $W$, since $W$ is quadratic and the sum of all angles
at every vertex of the tesselation is $2 \pi$.

The genus of the Wicks form $W$ is $m-1$ and the hyperbolic plane
tesselated by $D$ is a periodic plane under the action of genus 
$g=m-1$ surface group, the graph of $W$ is on the fig.5.

The Theorem is proved.

\begin{picture}(300,150)
\put(50,80){\circle*4}
\put(150,80){\circle*4}
\put(250,80){\circle*4}
\put(100,80){\oval(100,100)}
\put(100,80){\oval(100,50)[t]}
\put(200,80){\oval(100,100)}
\put(200,80){\oval(100,50)[t]}
\put(95,30){\vector(1,0){10}}\put(100,35){$a_m$}
\put(95,130){\vector(1,0){10}}\put(100,135){$a_1$}
\put(95,105){\vector(1,0){10}}\put(100,110){$a_2$}
\put(195,30){\vector(1,0){10}}\put(200,35){$b_m$}
\put(195,130){\vector(1,0){10}}\put(200,135){$b_1$}
\put(195,105){\vector(1,0){10}}\put(200,110){$b_2$}
\put(50,80){\vector(1,0){55}}\put(105,80){\line(1,0){45}}
\put(150,80){\vector(1,0){55}}\put(205,80){\line(1,0){45}}
\multiput(102,65)(0,-5){3}{\circle*{2}}
\multiput(202,65)(0,-5){3}{\circle*{2}}
\put(130,5){\bf{Figure 5.}}
\end{picture}

\end{document}